\documentclass[a4paper,11pt]{article}
\usepackage{amsmath,amssymb}
\usepackage{geometry}
\usepackage{graphicx}
\usepackage[colorlinks=true, urlcolor=blue, citecolor=black, linkcolor=black]{hyperref}
\usepackage{listings,xcolor}
\allowdisplaybreaks
\newcommand{\R}{\mathbb{R}}
\newcommand{\mv}[1]{\mathbf{#1}}
\newcommand{\e}[1]{\mv{e}_{#1}}
\newcommand{\rev}[1]{\widetilde{#1}}

\title{Comparison of Methods for
  Rotating a Point in~\texorpdfstring{\(\R^3\)}{R^3}: \\
  From Vector Algebra to Geometric Algebra}
\author{Tom Verhoeff\\
  Dept.\ Math.\ \& CS, Eindhoven University of Technology, The Netherlands\\
  \href{mailto:T.Verhoeff@tue.nl}{\ttfamily T.Verhoeff@tue.nl}}
\date{11~Feb 2025; updated 05~Apr 2025 and revised 20~Jun 2026}

\begin{document}
\maketitle

\begin{abstract}
\noindent
This article starts by presenting and comparing four methods
for computing
the rotation of a point about an axis by an angle in~\(\R^3\).
We illustrate these methods by computing, by hand, the rotation of
point~\(P=(1,0,1)^T\)
about axis~\(\mv{a}=(1,1,1)^T\)
by angle~\(\theta=60^\circ\) (following the right-hand rule).
The four methods considered are:
\begin{enumerate}
  \item An \emph{ad hoc} geometric method
    exploiting a symmetry in the situation.
  \item A projection method that sets up a new coordinate system
    using the dot and cross products.
  \item A matrix method which rotates the standard basis and
    uses matrix--vector multiplication.
  \item A Geometric Algebra method that represents the rotation
    as a double reflection via a rotor.
    We provide a brief introduction to Geometric Algebra.
\end{enumerate}
Subsequently,
we also address rotations in other dimensions and
explain how these can be handled with Geometric Algebra to provide
deeper insights.

\end{abstract}

\section{Introduction}
Rotations in three-dimensional space occur in many areas of
science and engineering.
This article originated as an attempt to understand more deeply
how Geometric Algebra (GA) handles rotations and what the pedagogical
challenges might be when introducing GA early
in the curriculum.
Working through a single concrete example and
comparing it with more familiar approaches
turned out to be surprisingly instructive.
The resulting comparison may also be useful for
members of the Pre-University Geometric Algebra group on LinkedIn\footnote{%
\url{https://www.linkedin.com/groups/8278281/}},
for mathematically interested secondary-school students, and
for teachers looking for a first encounter with rotors and GA.

The present version substantially revises and extends an earlier arXiv note.
Besides correcting several computational and conceptual errors in that version,
it adds a broader discussion of how the geometric-algebra treatment
of rotations relates to complex numbers in two dimensions,
quaternions in three dimensions, and
rotations in higher-dimensional spaces.

\begin{figure}[hbt]
\centering
{\includegraphics[trim=0 1.5cm 0 3.5cm,clip,height=6cm]{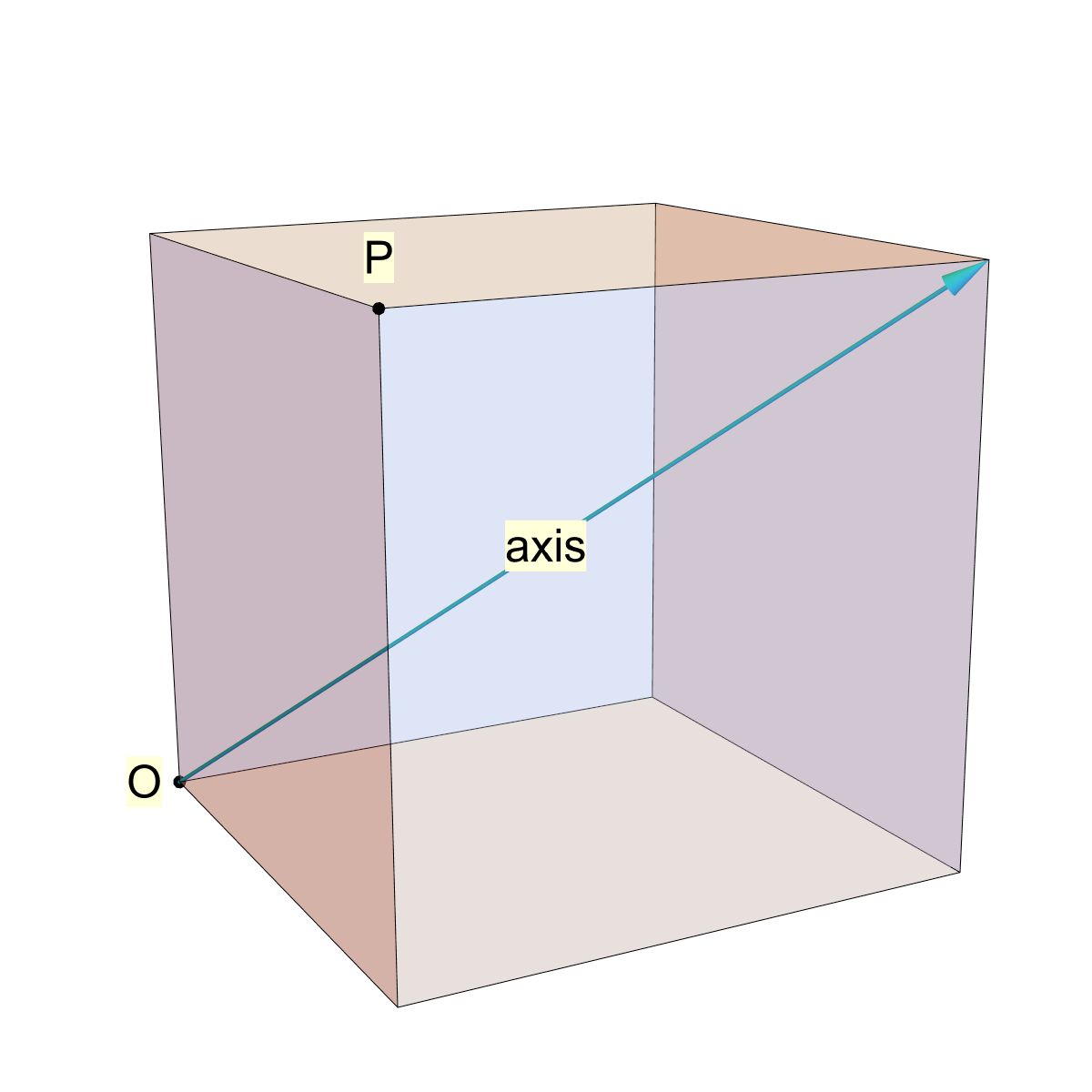}}
\caption{Rotation problem: Rotate given point~\(P\) about given axis \(\mv{a}\) by given angle~\(\theta\).}
\label{fig:set-up}
\end{figure}

We examine four methods for performing a rotation in~\(\R^3\)
``by hand''.
We shall use as our running example (see Fig.~\ref{fig:set-up})
the rotation of point (position vector)
\[P=(1,0,1)^T\]
about the line in the direction
\[\mv{a}=(1,1,1)^T\]
by an angle \[\theta=60^\circ\] (using the right-hand rule).
The first three methods require only elementary vector algebra and
linear algebra in~\(\R^3\),
in particular, the dot product and cross product,
but all definitions will be given.
The fourth method introduces the basic ideas of Geometric Algebra,
sometimes also called Clifford Algebra.

We restrict ourselves to rotation axes passing through the origin\footnote{%
For an axis not through the origin, translate point and axis
such that the axis passes through the origin,
rotate, then translate back.}
and we identify\footnote{%
In \emph{Conformal} Geometric Algebra (CGA),
one can represent points directly,
but that machinery is not needed here.}
point~\(P\) with its position vector~\(\overrightarrow{OP}\).
Thus, rotating the point means rotating its position vector.

\noindent
Although the particular configuration allows some shortcuts,
our aim is to compare the underlying concepts and
the conceptual load involved in each approach.

\section{Method 1: \emph{Ad Hoc} Approach Exploiting Symmetry}
\begin{figure}[hbt]
\centering
{\includegraphics[trim=0 3.5cm 0 2.5cm,clip,height=6cm]{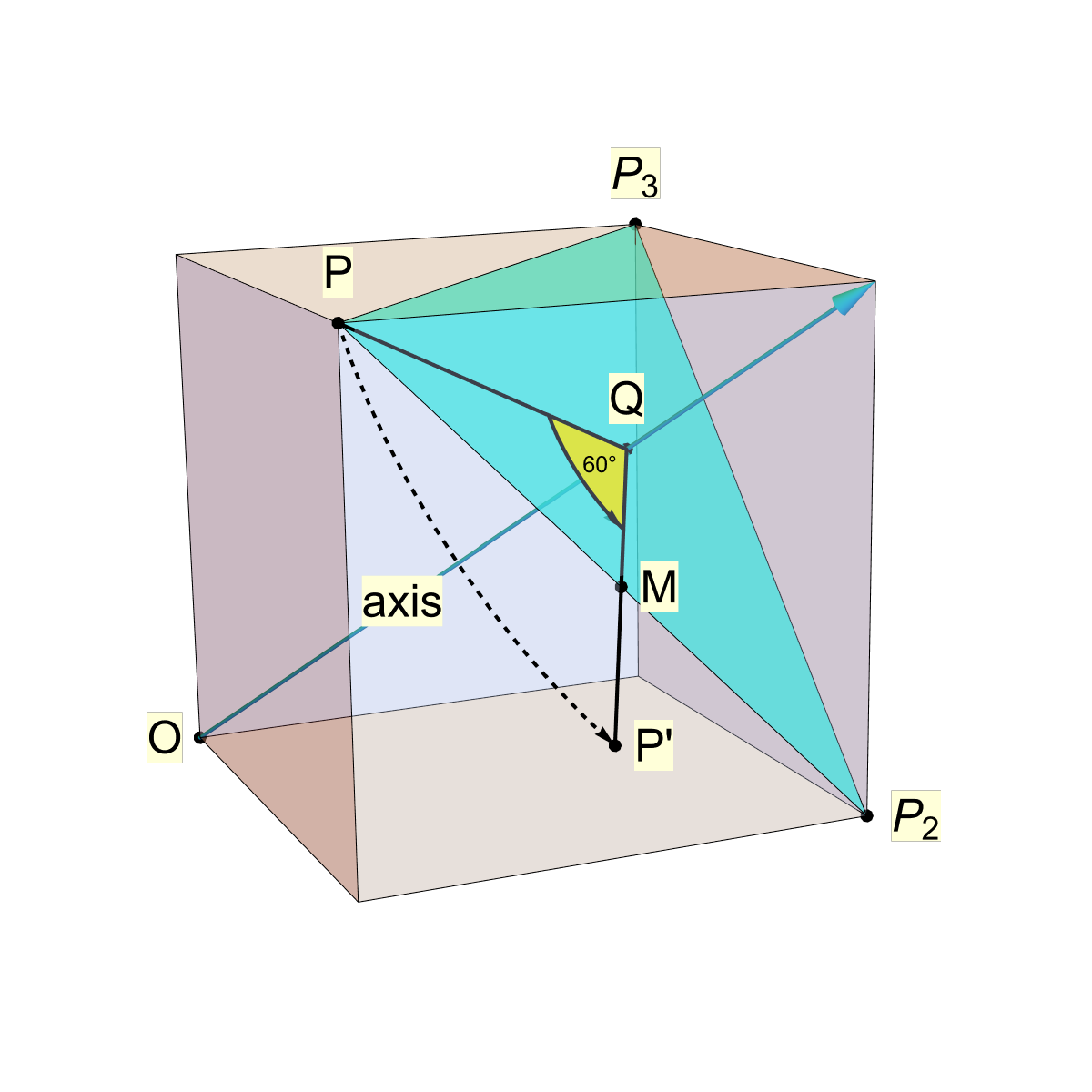}}
\caption{Method~1: \emph{Ad hoc} approach exploiting symmetry.}
\label{fig:ad-hoc}
\end{figure}

Method 1 is included mainly because this particular example
admits an elegant geometric shortcut.
Observe (see Fig.~\ref{fig:ad-hoc})
that the plane perpendicular to the axis \(\mv{a}=(1,1,1)^T\) and
containing \(P=(1,0,1)^T\)
is given by
\[
x+y+z=2,
\]
since
\[
1+0+1=2.
\]

\noindent
It turns out that the points
\[
P=(1,0,1)^T,\quad P_2=(1,1,0)^T,\quad P_3=(0,1,1)^T
\]
all lie in this plane and, moreover, one may verify that
\[
\|P-P_2\|=\|P_2-P_3\|=\|P_3-P\|=\sqrt{2}.
\]
Thus these three points form an \emph{equilateral triangle}. Its centroid (which, for an equilateral triangle, is also the circumcenter) is
\[
Q = \left(\frac{1+1+0}{3},\,\frac{0+1+1}{3},\,\frac{1+0+1}{3}\right)^T = 
  \left(\tfrac{2}{3},\tfrac{2}{3},\tfrac{2}{3}\right)^T .
\]
Since \(Q\) lies on axis~\(\mv{a}\), this point~\(Q\)
is also the projection of~\(P\) onto~\(\mv{a}\)
(also see Method~2 below).
Since the vertices of the triangle are separated by central angles of \(120^\circ\),
the symmetric configuration can be exploited to construct
a \(60^\circ\) rotation about the axis geometrically.
Namely, \(P'\) ends up on the line from~\(Q\) to the midpoint~\(M\) of~\(P P_2\),
such that \(\|P'-Q\|=\|P-Q\|=\sqrt{\frac{2}{3}}\).
A short calculation shows that
\[
M = \left(\frac{1+1}{2},\, \frac{0+1}{2},\, \frac{1+0}{2}\right)^T =
  \left(1,\tfrac{1}{2},\tfrac{1}{2}\right)^T
\]
and hence
\[
P'= Q + (M - Q)\frac{\|P-Q\|}{\|M-Q\|} =
  \left(\tfrac{4}{3},\tfrac{1}{3},\tfrac{1}{3}\right)^T.
\]

\noindent
For further discussion of symmetric configurations in space,
see, for example,~\cite{coxeter}.

\section{Method 2: Projection and a New Coordinate System}
\begin{figure}[hbt]
\centering
{\includegraphics[trim=0 2.9cm 0 3.7cm,clip,height=6cm]{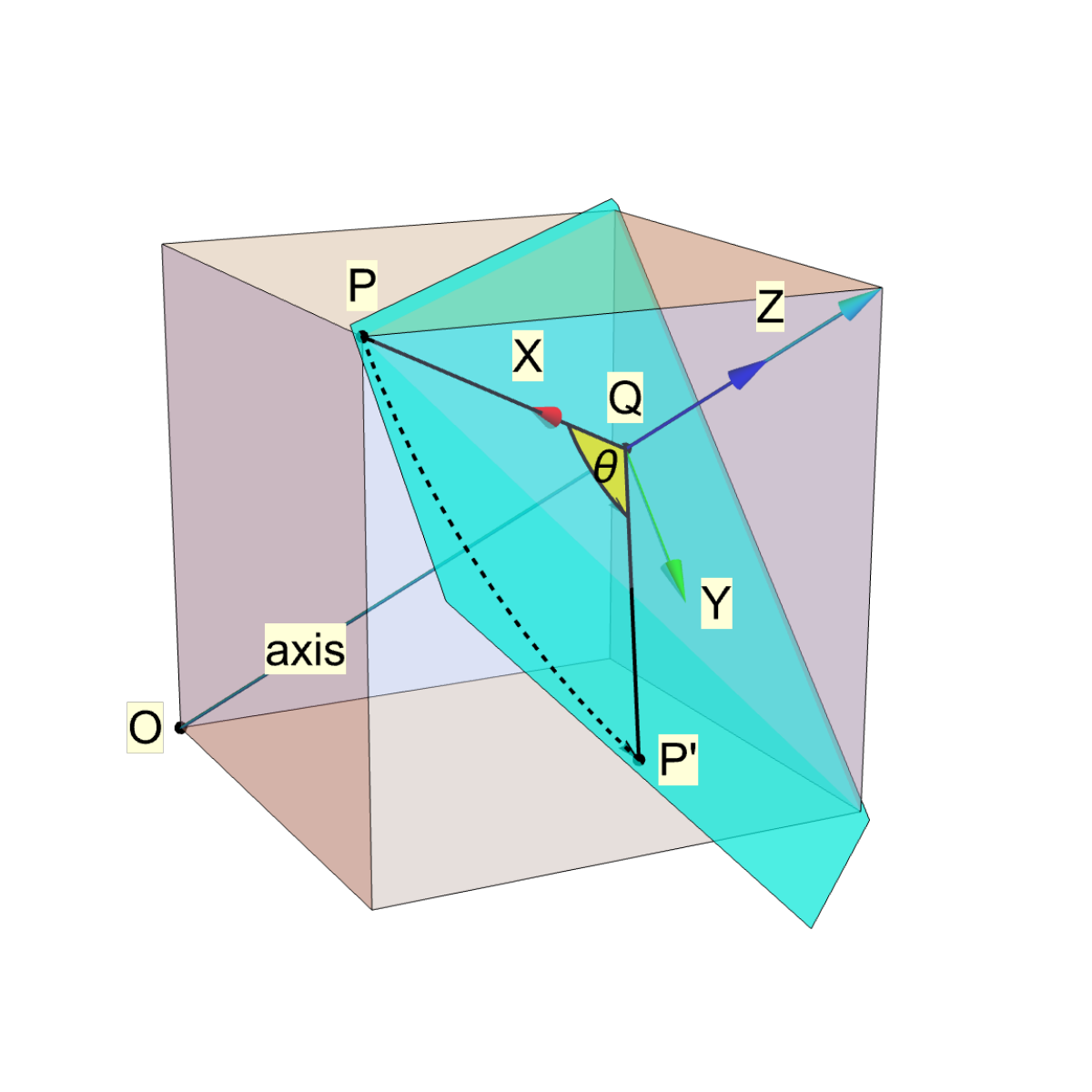}}
\caption{Method~2: Projection and new coordinate system
  (the latter is not to scale).}
\label{fig:projection}
\end{figure}

This method works for arbitrary rotation axes, rotation angles, and points.
See Fig.~\ref{fig:projection}.

\subsection{Projecting \texorpdfstring{\(P\)}{P} onto~\texorpdfstring{\(a\)}{a}}
The projection of \(P\) onto \(\mv{a}\) is given by
\[
Q=\frac{P\cdot \mv{a}}{\mv{a}\cdot \mv{a}}\,\mv{a}.
\]
Recall that the \emph{dot product} for vectors \(\mv{v}=(v_x,v_y,v_z)\) and
\(\mv{w}=(w_x,w_y,w_z)\) is
\[
\mv{v}\cdot \mv{w} = v_x w_x + v_y w_y + v_z w_z,
\]
and geometrically
\[
\mv{v}\cdot \mv{w}=\|\mv{v}\|\|\mv{w}\|\cos\psi,
\]
where \(\psi\) is the angle between~\(\mv{v}\) and~\(\mv{w}\).
For
\[
P=(1,0,1)^T \quad \text{and} \quad \mv{a}=(1,1,1)^T
\]
we have
\[
P\cdot \mv{a} = 1+0+1=2,\quad \mv{a}\cdot\mv{a} = 1+1+1=3,
\]
so that
\[
Q=\tfrac{2}{3}(1,1,1)^T = \left(\tfrac{2}{3},\tfrac{2}{3},\tfrac{2}{3}\right)^T .
\]

\subsection{Setting up a new coordinate system in the plane}
The vector
\[
P-Q=\left(1-\tfrac{2}{3},\,0-\tfrac{2}{3},\,1-\tfrac{2}{3}\right)^T
=\left(\tfrac{1}{3},-\tfrac{2}{3},\tfrac{1}{3}\right)^T
\]
lies in a plane perpendicular to axis~\(\mv{a}\).
We choose this direction as our new \(X\)--axis:
\[
  \mv{x} = \frac{P - Q}{\|P - Q\|} =
    \frac{\left(\tfrac{1}{3},-\tfrac{2}{3},\tfrac{1}{3}\right)^T}{\sqrt{6}/3}
    =\sqrt{\tfrac{3}{2}}\left(\tfrac{1}{3},-\tfrac{2}{3},\tfrac{1}{3}\right)^T .
\]
Next, we take the new \(Z\)--axis to be in the direction of~\(\mv{a}\), i.e.,
\[
\mv{z} = \frac{\mv{a}}{\|\mv{a}\|}=\tfrac{1}{\sqrt{3}}(1,1,1)^T,
\]
and define the \(Y\)--axis via the \emph{cross product}:
\[
\mv{y} = \mv{z}\times\mv{x}.
\]
For vectors \(\mv{v}= (v_x,v_y,v_z)\) and \(\mv{w}=(w_x,w_y,w_z)\)
the cross product is defined as
\[
\mv{v} \times \mv{w} = (v_y w_z - v_z w_y,\, v_z w_x - v_x w_z,\, v_x w_y - v_y w_x),
\]
and geometrically by
\(\mv{v} \times \mv{w}\) is perpendicular to both~\(\mv{v}\) and~\(\mv{w}\),
using the right-hand rule for direction,
and \(\|\mv{v} \times \mv{w}\|=\|\mv{v}\|\|\mv{w}\|\sin\phi\) where
\(\phi\) is the angle between~\(\mv{v}\) and~\(\mv{w}\).

A short calculation shows
\[
\mv{y}
=\mv{z} \times \mv{x}
=\tfrac{1}{\sqrt{3}}(1,1,1)^T\times
  \sqrt{\tfrac{3}{2}}\left(\tfrac{1}{3},-\tfrac{2}{3},\tfrac{1}{3}\right)^T
=\tfrac{1}{\sqrt{2}}\left(1,0,-1\right)^T.
\]

\subsection{Rotating in the new coordinate system}
Express the rotated vector (with respect to the center \(Q\)) as
\[
\frac{P'-Q}{\|P-Q\|}=\cos60^\circ\,\mv{x} + \sin60^\circ\,\mv{y} .
\]
With \(\cos60^\circ=\frac{1}{2}\) and \(\sin60^\circ=\frac{\sqrt{3}}{2}\),
we get
\[
P'-Q=\tfrac{\sqrt{6}}{3}\left(\tfrac{1}{2}\sqrt{\tfrac{3}{2}}\left(\tfrac{1}{3},-\tfrac{2}{3},\tfrac{1}{3}\right)^T
+\tfrac{\sqrt{3}}{2}\tfrac{1}{\sqrt{2}}\left(1,0,-1\right)^T\right) .
\]
A short computation gives
\[
P'-Q=\tfrac{1}{2}\left(\tfrac{1}{3}+1,\, -\tfrac{2}{3}+0,\, \tfrac{1}{3}-1\right)^T
=\left(\tfrac{2}{3},-\tfrac{1}{3},-\tfrac{1}{3}\right)^T .
\]
Thus,
\[
P'=Q+(P'-Q)
=\left(\tfrac{2}{3},\tfrac{2}{3},\tfrac{2}{3}\right)^T
+\left(\tfrac{2}{3},-\tfrac{1}{3},-\tfrac{1}{3}\right)^T
=\left(\tfrac{4}{3},\tfrac{1}{3},\tfrac{1}{3}\right)^T .
\]
This method is also easy to apply for other rotation angles~\(\theta\):
\[
P' = Q + \|P-Q\|(\mv{x}\cos\theta + \mv{y}\sin\theta) .
\]

\section{Method 3: The Matrix Approach}
\begin{figure}[hbt]
\centering
{\includegraphics[trim=0 0.5cm 0 1cm,clip,height=6cm]{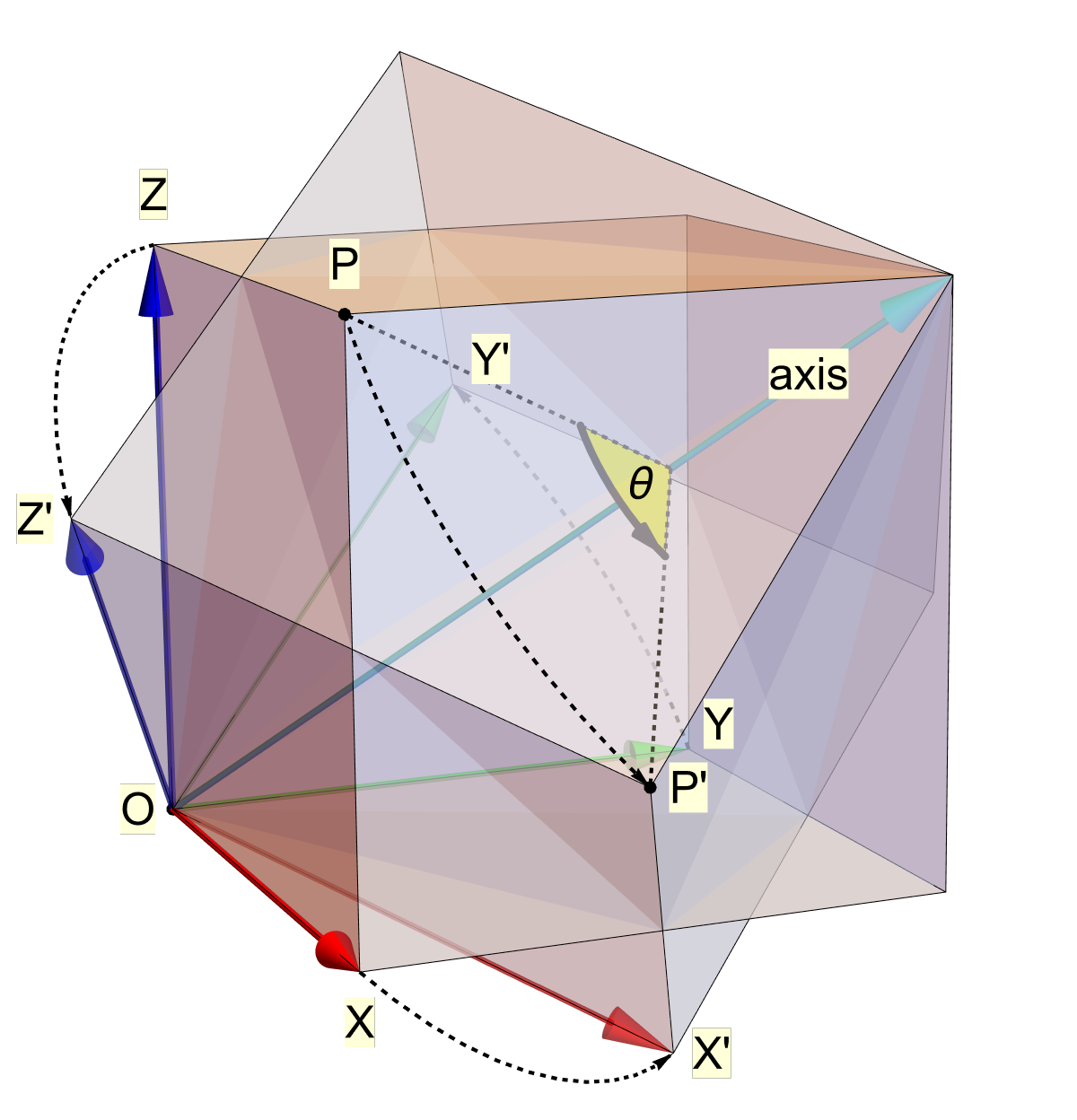}}
\caption{Method~3: Using a rotation matrix
  based on rotated basis vectors.}
\label{fig:rotation-matrix}
\end{figure}

When one needs to rotate various points~\(P=(x,y,z)^T\),
one can write each such point as a linear combination of
the basis vectors \(\e{1}, \e{2}, \e{3} = (1,0,0)^T, (0,1,0)^T, (0,0,1)^T\):
\[
P = x\, \e{1} + y\, \e{2} + z\, \e{3}
\]
Using Method~2, each basis vector can be rotated to yield
vectors \(\mv{e}'_1, \mv{e}'_2, \mv{e}'_3\)
(see Fig.~\ref{fig:rotation-matrix}).
Since rotation about an axis through the origin is a linear operation,
we now have
\[
P' = x\, \mv{e}'_1 + y\, \mv{e}'_2 + z\, \mv{e}'_3
\]
The latter is often written as a matrix-vector multiplication:
\[
P' = \begin{pmatrix}\mv{e}'_1\ \mv{e}'_2\ \mv{e}'_3\end{pmatrix}
  \begin{pmatrix} x\\ y\\ z\end{pmatrix}
\]
where the rotated basis vectors are written as columns.

It is cumbersome to rotate the three basis vectors individually by hand.
Any rotation about a fixed axis through the origin can be written as a
matrix-vector multiplication.
One standard formula, viz.\ Rodrigues' rotation formula in matrix form,
for the rotation matrix~\(R_{\text{mat}}\) about a unit vector 
\[
\mv{u} = (u_x,u_y,u_z)
\]
by an angle~\(\theta\) is given by
\[
\begin{aligned}
R_{\text{mat}} & = I + (\sin\theta)\,[\mv{u}]_\times + (1-\cos\theta)\,[\mv{u}]_\times^2 \\
  & = \cos\theta\,I+(1-\cos\theta)\,\mv{u}\mv{u}^T+(\sin\theta)\,[\mv{u}]_\times ,
\end{aligned}
\]
where \(I\) is the identity matrix and
\[
[\mv{u}]_\times = \begin{pmatrix}
0 & -u_z & u_y\\[1mm]
u_z & 0 & -u_x\\[1mm]
-u_y & u_x & 0
\end{pmatrix}
\]
is the skew-symmetric matrix corresponding to the cross product,
i.e., \(\mv{u}\times \mv{v} = [\mv{u}]_\times\,\mv{v}\).
The matrix \(\mv{u}\mv{u}^T=I+[\mv{u}]_\times^2\) is the \emph{outer product}.
See, e.g.,~\cite{dai}
for further details on the derivation of this formula.

For our case, since \(\mv{a}=(1,1,1)^T\) we have
\[
\mv{u} = \tfrac{1}{\sqrt{3}}(1,1,1)^T,
\]
and with \(\theta=60^\circ\) (so that \(\cos60^\circ=\frac{1}{2}\) and \(\sin60^\circ=\frac{1}{2}\sqrt{3}\)) the entries of the rotation matrix
can be computed to yield
\[
R_{\text{mat}} =
\tfrac{1}{2}\begin{pmatrix}1 & 0 & 0 \\ 0 & 1 & 0 \\ 0 & 0 & 1 \end{pmatrix} +
\tfrac{1}{6}\begin{pmatrix}1 & 1 & 1 \\ 1 & 1 & 1 \\ 1 & 1 & 1 \end{pmatrix} +
\tfrac{1}{2}\begin{pmatrix}\phantom{-}0 & -1 & \phantom{-}1 \\ \phantom{-}1 & \phantom{-}0 & -1 \\ -1 & \phantom{-}1 & \phantom{-}0 \end{pmatrix}
=
\begin{pmatrix}
\phantom{-}\frac{2}{3} & -\frac{1}{3} & \phantom{-}\frac{2}{3}\\[1mm]
\phantom{-}\frac{2}{3} & \phantom{-}\frac{2}{3} & -\frac{1}{3}\\[1mm]
-\frac{1}{3} & \phantom{-}\frac{2}{3} & \phantom{-}\frac{2}{3}
\end{pmatrix}.
\]
Then the rotated point is
\[
P'=R_{\text{mat}}\,P
=\begin{pmatrix}
\frac{2}{3}\cdot1+(-\frac{1}{3})\cdot0+\frac{2}{3}\cdot1\\[1mm]
\frac{2}{3}\cdot1+\frac{2}{3}\cdot0+(-\frac{1}{3})\cdot1\\[1mm]
-\frac{1}{3}\cdot1+\frac{2}{3}\cdot0+\frac{2}{3}\cdot1
\end{pmatrix}
=\begin{pmatrix}
\frac{4}{3}\\[1mm]
\frac{1}{3}\\[1mm]
\frac{1}{3}
\end{pmatrix}.
\]
Note that a rotation matrix operates on vectors and
that rotation operations can be composed,
without having to apply them to a vector,
via \emph{matrix-matrix multiplication} of their rotation matrices.
Storing a rotation matrix requires nine real numbers.
Since there are only three rotational degrees of freedom in~3D,
those nine numbers must satisfy orthogonality constraints,
which is expensive to maintain numerically under composition.

\section{Method 4: The Geometric Algebra Approach}
Whereas Methods~2 and~3 describe rotations using coordinates and matrices,
Geometric Algebra treats rotations as compositions of reflections.
We first provide the relevant definitions.
In (Euclidean) Geometric Algebra for~\(\R^3\),
\emph{multivectors}
are obtained as linear combinations of
products of the basis vectors \(\e{1}, \e{2}, \e{3}\).
The empty product corresponds to a scalar:
\[
  c_0 +
  c_1\, \e{1} + c_2\, \e{2} + c_3\, \e{3} +
  c_4\, \e{12} + c_5\, \e{23} + c_6\, \e{31} +
  c_7\, \e{123}
\]
where \(\e{ij}\) is short for the product~\(\e{i} \e{j}\).
The number of basis vectors in such a product is called its \emph{grade}.
Thus scalars have grade~\(0\), ordinary vectors have grade~\(1\),
bivectors have grade~\(2\), and trivectors have grade~\(3\).
A general multivector may contain terms of several grades.
A multivector containing only even grades,
such as grades~\(0\) and~\(2\), is called \emph{even}.

Multivectors represent \emph{geometric objects}.
The scalar component~\(c_0\) can be viewed as a sized and oriented
0-dimensional subspace and
the \(c_i\, \e{i}\) components are sized and oriented 1-dimensional subspaces
(size is vector length).
Similarly,
the \(c_k\, \e{ij}\) components are sized and oriented (but unshaped)
2-dimensional subspaces
(size is area;
oriented also means that the plane segment has distinct front and back sides).
Finally, \(c_7\, \e{123}\) is a sized and oriented 3-dimensional subspace
(size is volume; there is only one orientation modulo its sign;
hence,
\(\e{123}=\mv{I}\) is also known as the space's \emph{pseudoscalar}).
Any vector corresponds to a linear combination of
the three \(\e{i}\) basis vectors.
Similarly, any oriented plane segment in~\(\R^3\) corresponds to
a linear combination of the three basis \emph{bivectors}~\(\e{ij}\).

Products of basis vectors satisfy the following relations:
\begin{align*}
  \e{i} \e{i} & = 1 \\
  \e{i} \e{j} & = - \e{j} \e{i} \text{ if } i \ne j
\end{align*}
Note that addition and multiplication are associative,
that addition is commutative, but that multiplication is not.
Moreover, multiplication distributes over addition.
Multiplication of multivectors is known as the \emph{geometric product}.
For example, we have
\[
  \left(\e{12}\right)^2 = \e{1}\e{2}\e{1}\e{2} = - \e{1}\e{1}\e{2}\e{2} = -1
\]
and thus the bivector \(\e{12}\)
behaves like the complex unit~\(i\) with~\(i^2=-1\).
Also see Appendix~\ref{sec:GA-details} for a few more examples.

For our rotation problem, we can write
\begin{align*}
\mv{a} & = \e{1} + \e{2} + \e{3} \\
P & = \e{1} + \e{3}
\end{align*}
A 3D rotation
can be accomplished by two successive \emph{reflections}.
Reflection of vector~\(\mv{v}\) in the plane
perpendicular to unit vector~\(\mv{b}\)
can be expressed using the geometric product as
(see Appendix~\ref{sec:GA-reflection} for an explanation)
\[
  v \longmapsto - \mv{b} \mv{v} \mv{b} .
\]
This is known as a \emph{sandwich product},
where \(\mv{v}\) is sandwiched between two~\(\mv{b}\)'s.
Hence,
reflecting first in the plane with normal~\(\mv{b}\)
and then in the plane with normal~\(\mv{c}\) gives
\[
  v \longmapsto -\mv{c}(-\mv{b} \mv{v} \mv{b}) \mv{c}
  = \mv{c} \mv{b} \mv{v} \mv{b} \mv{c} .
\]
The geometry of this transformation is classical:
the composition of two reflections in planes through the origin is a rotation.
Its axis is the line of intersection of the two planes,
and its rotation angle is twice the angle between the planes.

Conversely,
suppose we want a rotation through angle~\(\theta\)
about a given axis through the origin.
Then we could choose two planes containing
that axis and making angle~\(\theta/2\),
and use their corresponding unit normal vectors~\(\mv{b}\) and~\(\mv{c}\).

The product
\[
        R = \mv{c}\mv{b}
\]
called a \emph{rotor},
then encodes the rotation, and the rotated vector is
\[
        \mv{v}' = R\mv{v}\rev{R} ,
\]
where \(\rev{R}=\mv{b}\mv{c}\) denotes reversal.

For the present computation, however,
we do not actually need to construct the two reflecting planes.
The useful fact, explained in Appendix~\ref{sec:GA-rotor},
is that in three dimensions every rotation about
the origin can be represented by a rotor of the form
\[
        R
        =
        \cos\tfrac{\theta}{2}
        -
        \mv{B}\sin\tfrac{\theta}{2}
        =
        e^{-\mv{B}\theta/2},
\]
where \(\mv{B}\) is the unit bivector representing the
\emph{plane of rotation}.
This is the same kind of rotor as the product \(\mv{c}\mv{b}\) above,
but written directly in terms of the rotation angle and
the rotation plane.

Thus, the remaining task is to find the appropriate unit bivector~\(\mv{B}\).
The given rotation axis has unit direction
\[
        \mv{n}=\tfrac{1}{\sqrt3}(\e1+\e2+\e3).
\]
In three dimensions, the pseudoscalar \(\mv{I}=\e1\e2\e3\)
can be used to convert an axis direction
into the perpendicular rotation plane
(cf.\ Appendix~\ref{sec:GA-perpendicular-bivector}).
With our orientation convention we take
\[
        \mv{B}=\mv{I}\mv{n}
        =
        \tfrac{1}{\sqrt3}\e{123}(\e1+\e2+\e3)
        =
        \tfrac{1}{\sqrt3}(\e{23}+\e{31}+\e{12}) .
\]
Note that \(\e{123}(\e{1} + \e{2} + \e{3}) =
    \e{1231} + \e{1232} + \e{1233} = \e{1123} - \e{1223} + \e{1233} =
    \e{23} + \e{31} + \e{12}\).
Therefore the rotor for the desired rotation through angle~\(\theta\) is
\[
        R
        =
        \cos\tfrac{\theta}{2}
        -
        \tfrac{1}{\sqrt3}(\e{23}+\e{31}+\e{12})
        \sin\tfrac{\theta}{2} .
\]
\begin{figure}[hbt]
\centering
{\includegraphics[trim=0 3.5cm 0 4cm,clip,height=5cm]{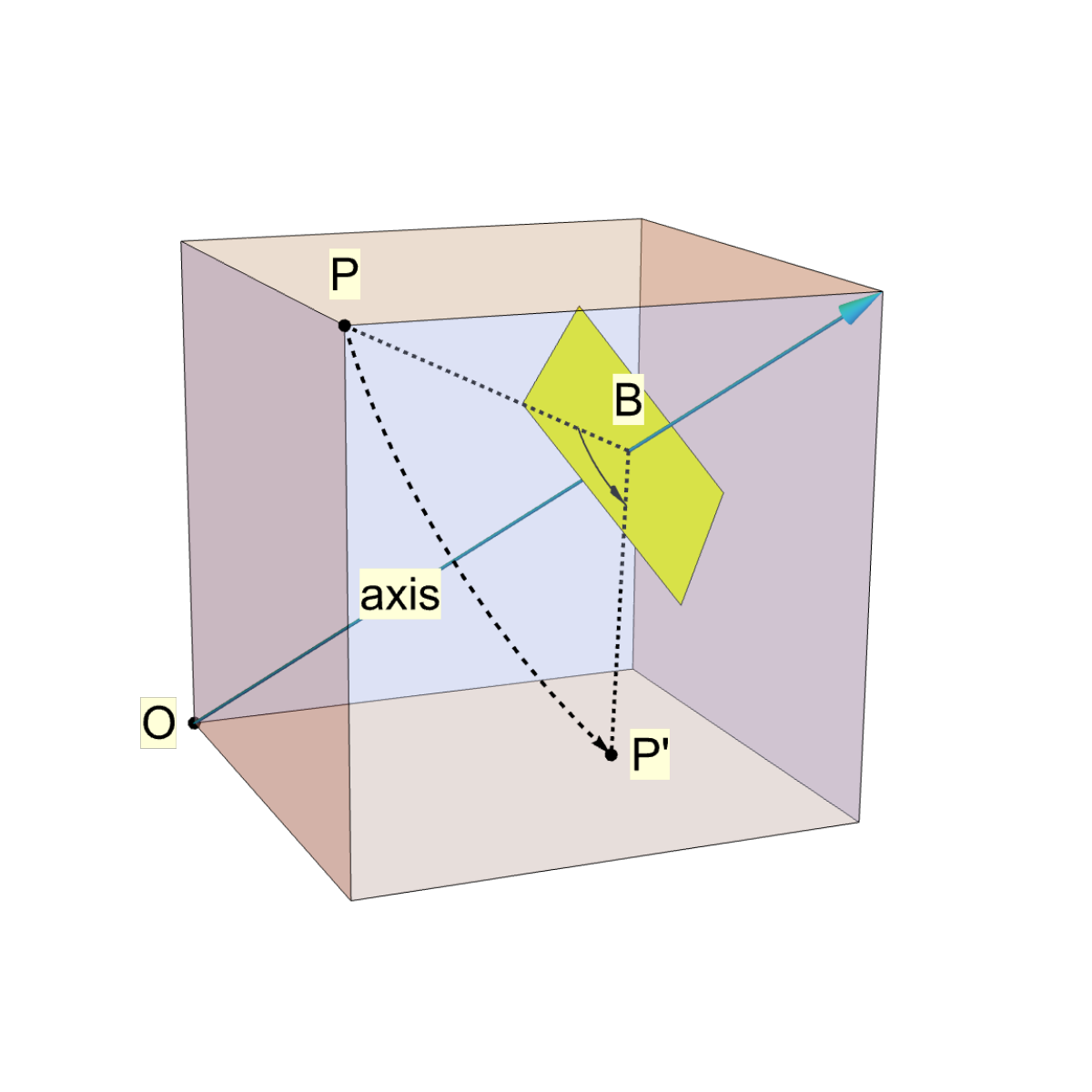}}
\caption{Method~4: Geometric Algebra
  (unit bivector~\(\mv{B}\) is not to scale).}
\label{fig:GA}
\end{figure}

\noindent
With \(\theta=60^\circ\), so that \(\theta/2=30^\circ\),
the rotor becomes
\[
R=\cos30^\circ-\mv{B}\sin30^\circ
=\tfrac{\sqrt{3}}{2}-\tfrac{1}{2}\,\mv{B}
\quad\text{and}\quad \rev{R}=\tfrac{\sqrt{3}}{2}+\tfrac{1}{2}\,\mv{B}.
\]
Substituting the rotor into the sandwich product and
expanding using the multiplication rules yields
(see Appendix~\ref{sec:GA-calculation} for the full algebraic calculation)
\[
P' = R\,P\,\rev{R} =
  \tfrac{4}{3}\e{1} + \tfrac{1}{3} \e{2} + \tfrac{1}{3} \e{3} .
\]
Note that a rotor is an operator on vectors and
that rotation operations can be composed,
without having to apply them to a vector,
via the geometric product of the rotors.
Storing a rotor requires four real numbers.
Since there are only three rotational degrees of freedom in~3D,
there is a relationship between these numbers,
viz.\ that the magnitude of the rotor is~\(1\),
and this is easy to maintain under composition.
Note that 3D~rotors are equivalent to \emph{unit quaternions}
(cf.~\S\ref{sec:quaternions}).
For more on this approach see, e.g.,~\cite{hestenes} and~\cite{dorst}.

\section{Comparison and Discussion}
We have described four different methods for computing the rotation of
\[
P=(1,0,1)^T
\]
about the axis
\[
\mv{a}=(1,1,1)^T
\]
by angle~\(\theta=60^\circ\).
The first method uses a symmetry unique to the situation
(namely, the fact that \(P\) lies in an equilateral triangle
together with \((0,1,1)^T\) and \((1,1,0)^T\)).
The second method sets up a new coordinate system
using the projection of~\(P\) onto~\(a\)
and uses the dot and cross products (whose definitions we recalled)
to “rotate” the vector.
The third method derives the rotation matrix by “rotating” the standard basis;
this approach is particularly useful
if many vectors must be rotated about the same axis.
Finally, the fourth method, based on Geometric Algebra,
interprets the rotation as the product of two reflections
(or equivalently via the rotor \(R=\exp(-\mv{B}\theta/2)\)),
unifying the treatment of rotations and reflections.

Though each approach has its own conceptual overhead and computational flavor, all yield the same final result:
\[
P'=\left(\tfrac{4}{3},\tfrac{1}{3},\tfrac{1}{3}\right)^T.
\]
The choice of method depends on the context and the user’s familiarity with the underlying mathematical tools.

Some notes:
\begin{enumerate}
\item The \emph{ad hoc} approach works only in very specific situations
  (axis~\(\mv{a}\), angle~\(\theta\), vector~\(P\)).
  Conceptually, it requires the use of Cartesian coordinates, vectors,
  computing midpoints, vector length, and scaling.
\item The projection method is fully general,
  and varying angle~\(\theta\) is easy,
  whereas varying axis~\(\mv{a}\) or vector~\(P\) is more involved.
  Conceptually, it additionally requires
  the dot product and cross product of vectors,
  and trigonometry for the rotation.
\item The matrix method is also fully general,
  and it is very easy to vary vector~\(P\),
  whereas varying axis~\(\mv{a}\) or angle~\(\theta\)
  requires the computation of a new rotation matrix.
  
  Conceptually, this requires additionally
  the use of matrices and matrix-vector multiplication.
  The rotation matrix can be obtained by applying the projection method
  for each of the basis vectors,
  or by applying a (rather complicated) formula.
  
  Matrices are nice because rotations can be composed
  by multiplying their matrices, using matrix-matrix multiplication
  (something new to be learned).
  Maintaining orthogonality under repeated numerical composition
  may require re-or\-tho\-nor\-mal\-i\-za\-tion.
\item The Geometric Algebra (GA) method is fully general and
  works in spaces of arbitrary dimension
  (cf.~\S\ref{sec:other-dimensions}).
  The key insight is that, in general,
  a rotation can be captured by a rotor,
  which is a special kind of multivector.
  That rotor is used in an exponential and applied as a sandwich product.
  
  Conceptually, this requires basic understanding of GA concepts, such as
  multivectors, their algebra with addition and the geometric product.
  To understand why the method works, one also needs to study reflection
  and how products of vectors are composed of a scalar and a bivector,
  involving the angle between the vectors.
  Understanding the dot product and wedge product can be useful,
  but the cross product and matrices are not needed.
  The geometric product allows linear operators like reflection and rotation
  to be expressed within the algebra.

  The manual calculation with a rotor in a sandwich product is
  lengthy but elementary, and easy to automate.
  
  Rotors are nice because rotations can be composed
  by multiplying their rotors using the geometric product
  (which is already familiar).
  Numerically, the composition of rotors is easy to re-normalize.
  For the special case of~\(\R^3\),
  the GA method is equivalent to the use of quaternions.
\end{enumerate}

The four methods considered above all solve
the same three-dimensional rotation problem.
The geometric-algebra approach differs from the other methods
in one important respect.
While the first three methods are tied rather closely to
the special features of three-dimensional space,
the rotor formalism extends naturally to arbitrary dimensions,
where not every rotation has an axis.
The next section briefly explores this broader perspective and
its connections to complex numbers, quaternions, and
higher-dimensional rotations.

\section{Rotations in Other Dimensions:
The Geometric-Algebra Perspective}
\label{sec:other-dimensions}
This section is optional and provides background for readers
who want the broader payoff of geometric algebra.

So far, we dealt with rotations in three-dimensional Euclidean space
(with the rotation axis passing through the origin).
In that setting,
every rotation is a so-called \emph{simple rotation}:
in each plane perpendicular to the axis,
it acts as the same two-dimensional rotation;
thus, it does not move points on the axis.
Geometric algebra makes this situation especially transparent,
because such a rotation can be obtained as the composition of two reflections.

\subsection{Simple rotations from two reflections
in \texorpdfstring{\(n\)}{n}~dimensions}
The rotor construction used in Method~4 is not specific to three dimensions.
We therefore revisit it in a dimension-independent setting.
Let \(\mv{b}\) and \(\mv{c}\) be unit vectors
in \(n\)-dimensional Euclidean space.
They may be viewed as the unit normals of two reflecting hyperplanes
through the origin.
Reflection of a vector~\(\mv{v}\) in the hyperplane
with unit normal~\(\mv{b}\) is given by
\[
        \mv{v} \longmapsto -\mv{b}\mv{v}\mv{b} .
\]
Reflecting first in the hyperplane with normal~\(\mv{b}\),
and then in the hyperplane with normal~\(\mv{c}\), gives
\[
        \mv{v} \longmapsto \mv{c} \mv{b} \, \mv{v} \, \mv{b} \mv{c} .
\]
Thus the rotor is
\[
        R = \mv{c}\mv{b},
        \qquad
        \rev{R} = \mv{b}\mv{c},
\]
and the rotated vector is
\[
        \mv{v}' = R \mv{v} \rev{R} .
\]

Since \(\mv{b}\) and \(\mv{c}\) are vectors,
their geometric product decomposes as
\[
        \mv{c}\mv{b} = \mv{c}\cdot \mv{b} + \mv{c}\wedge \mv{b}
        = \mv{b}\cdot \mv{c} - \mv{b}\wedge \mv{c} .
\]
Thus,
this particular rotor consists of a scalar part \(\mv{b}\cdot \mv{c}\) and
a bivector part \(-\mv{b}\wedge \mv{c}\).
The bivector \(\mv{b}\wedge \mv{c}\) is called \emph{simple},
because it is explicitly the outer product of two vectors.
Geometrically, a simple bivector spans a single oriented plane.
(Note that e.g.\ bivector \(\e{12} + \e{34}\) is not simple
and combines two independent planes.)
If \(\mv{b}\) and \(\mv{c}\) are not parallel, we may write
\[
        \mv{B} = \frac{\mv{b}\wedge \mv{c}}{\|\mv{b}\wedge \mv{c}\|}
\]
for the corresponding unit bivector.
If the angle from~\(\mv{b}\) to~\(\mv{c}\) is~\(\varphi\), then
\[
        \mv{c}\mv{b} = \cos\varphi - \mv{B}\sin\varphi .
\]
The resulting rotation angle is~\(2\varphi\),
because a composition of two reflections in hyperplanes
meeting at angle~\(\varphi\) is a rotation through twice that angle.

Equivalently, if \(\mv{B}\) is already chosen as a simple unit bivector
representing the oriented plane of rotation,
then a rotation through angle~\(\theta\) can be written as rotor
\[
        R = e^{-\frac{\theta}{2}\mv{B}}
          = \cos\tfrac{\theta}{2}
            - \mv{B}\sin\tfrac{\theta}{2},
\]
and the action on vectors is again
\[
        \mv{v}' = R\mv{v}\rev{R} .
\]
The reverse is
\[
        \rev{R}
        =
        \cos\tfrac{\theta}{2}
        +
        \mv{B}\sin\tfrac{\theta}{2},
\]
and \(R\rev{R}=1\).

A rotor of the form
\[
        R = \cos\tfrac{\theta}{2}
            - \mv{B}\sin\tfrac{\theta}{2},
\]
with \(\mv{B}\) a unit \emph{simple} bivector,
represents a \emph{simple rotation}.
In \(n\)-dimensional space it rotates in the two-dimensional plane
represented by \(\mv{B}\) and
fixes the perpendicular \((n-2)\)-dimensional subspace.
In dimensions up to three, all bivectors are simple,
but from dimension four on,
there are also rotations that are not simple.

\subsection{Two dimensions and complex numbers}

In two dimensions, let
\[
        \mv{I} = \e{1}\e{2}
\]
be the unit pseudoscalar.
It satisfies
\[
        \mv{I}^2 = \e1\e2\e1\e2 = - \e1\e1\e2\e2 = -1 .
\]
The even subalgebra consists of elements
\[
        a + b \mv{I},
\]
and is therefore algebraically the same as the complex numbers.

It is important, however, not to confuse this with the representation of
vectors.
A vector in the plane is still represented as
\[
        \mv{v} = x \e{1} + y \e{2},
\]
not as \(x+y\mv{I}\).
The expression \(x+y\mv{I}\) is an even multivector,
whereas \(x \e{1}+y \e{2}\) is a vector.

A rotation through angle~\(\theta\) is represented by the rotor
\[
        R
        =
        e^{-\frac{\theta}{2}\mv{I}}
        =
        \cos\tfrac{\theta}{2}
        -
        \mv{I}\sin\tfrac{\theta}{2}
        \text{ with reverse }
        \rev{R}
        =
        \cos\tfrac{\theta}{2}
        +
        \mv{I}\sin\tfrac{\theta}{2}
        =
        e^{\frac{\theta}{2}\mv{I}}
        .
\]
The rotated vector is
\[
        \mv{v}' = R\mv{v}\rev{R} .
\]
Since \(\mv{I}\) anticommutes with every vector in the plane
(\(\mv{v}\mv{I} = - \mv{I}\mv{v}\)), we have
\[
        R\mv{v} = \mv{v}\rev{R} .
\]
Consequently the sandwich product simplifies to
\[
        \mv{v}' = \mv{v}\rev{R}^2
           =
        \mv{v} e^{\theta \mv{I}} .
\]
Thus, in two dimensions,
the usual multiplication by a unit complex number
with argument~$\theta$
appears naturally from the rotor formula.
The unit complex number is not a vector itself;
it is the operator that rotates vectors.

This also explains why two-dimensional rotations commute.
The even subalgebra is generated by~\(1\) and~\(\mv{I}\),
and is commutative.
Hence, all two-dimensional rotors commute with one another.

Although vectors and rotors play different roles in geometric algebra,
there is an intriguing similarity in two dimensions.
If we formally associate the vector
\[
  a\e1+b\e2
\]
with the complex number
\[
  a+b\mv{I},
\]
then multiplication by a rotor
\[
  c+d\mv{I}
\]
transforms the coefficients in exactly the same way:
\[
  (a\e1+b\e2)(c+d\mv{I})
=
  (ac-bd)\e1+(ad+bc)\e2,
\]
while
\[
  (a+b\mv{I})(c+d\mv{I})
  =
  (ac-bd)+(ad+bc)\mv{I}.
\]
This explains why complex numbers can be used successfully
to describe two-dimensional vectors and rotations.
From the geometric algebra point of view, however,
vectors and rotors remain fundamentally different objects:
vectors are linear combinations of~\(\e1\) and~\(\e2\),
whereas rotors are linear combinations of~\(1\) and~\(\mv{I}\).

\subsection{Three dimensions and quaternions}
\label{sec:quaternions}
In three dimensions, the even subalgebra is spanned by
\[
        1,\qquad \e{23},\qquad \e{31},\qquad \e{12}.
\]
These three bivectors all square to~\(-1\),
and their multiplication rules are,
up to a conventional choice of signs,
the multiplication rules for the quaternion units.
Thus, the even subalgebra of three-dimensional Euclidean
geometric algebra is isomorphic to the quaternion algebra.

Under this identification, a three-dimensional rotor
\[
        R
        =
        \cos\tfrac{\theta}{2}
        -
        \mv{B}\sin\tfrac{\theta}{2},
\]
where \(\mv{B}\) is a unit bivector, corresponds to a unit quaternion.
The quaternion ``imaginary vector part'' is, in geometric algebra terms,
a bivector.
In three dimensions, this distinction is often hidden by duality:
bivectors are dual to vectors.
For example, the bivector representing the plane
perpendicular to a unit axis vector~\(\mv{u}\) may be written using the
three-dimensional pseudoscalar \(\mv{I}_3=\e1\e2\e3\) as
\[
        \mv{B} = \mv{I}_3 \mv{u}
\]
up to the sign convention chosen for orientation.

Thus quaternions are not a different mechanism for rotations.
They are the three-dimensional even subalgebra of geometric algebra,
written in a dualized notation.
Geometric algebra keeps the geometric meaning explicit:
the generator of a rotation is a bivector,
representing the plane of rotation,
not primarily an \emph{axial} vector.

\subsection{Four dimensions: double rotations}

In dimensions higher than three, not every rotation is simple.
The product of two reflections always gives a simple rotation,
and hence a rotor with only a scalar part and a simple bivector part.
But according to the Cartan--Dieudonné theorem,
a general orthogonal transformation in \(n\)~dimensions
can be written as a composition of at most~\(n\) reflections.
Thus, 
in geometric algebra, a general rotor for a proper\footnote{%
An isometry is called \emph{proper} when it preserves handedness,
i.e., does not reflect.}
rotation is
an \emph{even} product of unit vectors:
\[
        R = a_{2k}\cdots a_2a_1 .
\]
Such a rotor is a unit even multivector,
but it need not consist only of a scalar and a simple bivector.

The first dimension where this distinction appears is four.
Consider the rotor
\[
        R
        =
        \big(
          \cos\tfrac{\alpha}{2}
          -
          \e{12}\sin\tfrac{\alpha}{2}
        \big)
        \big(
          \cos\tfrac{\beta}{2}
          -
          \e{34}\sin\tfrac{\beta}{2}
        \big),
\]
which is a composition of two simple rotations,
first with bivector~\(\e{12}\) over angle~\(\alpha\ne0\)
and second with bivector~\(\e{34}\) over~\(\beta\ne0\).
The two bivectors~\(\e{12}\) and~\(\e{34}\) represent orthogonal planes.
They commute:
\[
        \e{12}\e{34} = \e{34}\e{12} = \e{1234}.
\]
Expanding \(R\) gives
\[
        R
        =
        \cos\tfrac{\alpha}{2}\cos\tfrac{\beta}{2}
        \;-\;
        \e{12}\sin\tfrac{\alpha}{2}\cos\tfrac{\beta}{2}
        \;-\;
        \e{34}\cos\tfrac{\alpha}{2}\sin\tfrac{\beta}{2}
        \;+\;
        \e{1234}\sin\tfrac{\alpha}{2}\sin\tfrac{\beta}{2}.
\]
This rotor has grades~\(0\), \(2\), \emph{and}~\(4\).
So, a 4D~rotor need not merely be a scalar plus a simple bivector.

Its action on the basis vectors is easy to understand.
The first factor rotates vectors in the \(\e1\e2\)-plane
through angle~\(\alpha\),
while the second factor rotates vectors in the \(\e3\e4\)-plane
through angle~\(\beta\).
Thus, this is a simultaneous rotation in two orthogonal planes:
\[
\begin{aligned}
        \e1 &\longmapsto
        \e1\cos\alpha + \e2\sin\alpha, \\
        \e2 &\longmapsto
       -\e1\sin\alpha + \e2\cos\alpha, \\
        \e3 &\longmapsto
        \e3\cos\beta + \e4\sin\beta, \\
        \e4 &\longmapsto
       -\e3\sin\beta + \e4\cos\beta .
\end{aligned}
\]
That is
\[
  (x, y, z, w) \longmapsto
  (x\cos\alpha-y\sin\alpha,\, y\cos\alpha+x\sin\alpha,\,
   z\cos\beta-w\sin\beta,\, w\cos\beta+z\sin\beta)
\]
This is called a \emph{double rotation}.
In general, it fixes only the origin and moves every nonzero vector.
But vectors in the \(\e1\e2\)-plane (with \(z=w=0\)) remain in that plane,
i.e., the \(\e1\e2\)-plane is an \emph{invariant plane} of the rotation.
Also the \(\e3\e4\)-plane is an invariant plane (with \(x=y=0\)).

Let's investigate the \emph{orbits} of vectors under such rotations,
that is, the collection of vectors under rotation by~\(\alpha t\)
and~\(\beta t\) for~\(t\ge0\).
Consider the point with coordinates \((r_1, 0, r_2, 0)\),
which lies at distance~\(r_1\) from~\(\e{34}\) and
distance~\(r_2\) from~\(\e{12}\).
Its orbit is parameterized by
\[
  (r_1\cos(\alpha t),\, r_1\sin(\alpha t),\,
   r_2\cos(\beta t),\, r_2\sin(\beta t))
\]
The orbits of vectors in an invariant plane
(\(r_1=0\) or \(r_2=0\)) are ordinary circles.
In general, the orbits are helical curves on the 2D~surface
of a (Clifford) torus embedded in~\(\R^4\),
with radii~\(r_1\) and~\(r_2\).
When \(\alpha\) and~\(\beta\) have a rational ratio,
these orbits are closed curves (torus knots).
But when that ratio is irrational,
they lie dense on the torus surface and never close.

A special case occurs when
\[
        \alpha = \beta
        \qquad\text{or}\qquad
        \alpha = -\beta .
\]
The orbit of our example point simplifies to
\[
  (r_1, 0, r_2, 0)\cos(\alpha t) +
  (0, r_1, 0, \pm r_2)\sin(\alpha t)
\]
These are circles in the plane spanned by
\((r_1, 0, r_2, 0)\) and \((0, r_1, 0, \pm r_2)\).
In fact, all orbits are circles and
every vector is rotated through the same angle,
although in different two-dimensional planes depending on the vector.
Such rotations are called \emph{isoclinic rotations}.
Four-dimensional rotations
therefore have phenomena with no direct analogue in three dimensions.

\subsection{Summary}

The following distinctions are useful.
\begin{itemize}
\item
The product of two unit vectors,
\[
        \mv{c}\mv{b} = \mv{b}\cdot \mv{c} - \mv{b}\wedge \mv{c},
\]
is a rotor with scalar part and simple bivector part.
It represents a simple rotation.

\item
In three dimensions,
every rotation about the origin is simple.
Hence, every three-dimensional rotor can be written as
\(R=\cos(\theta/2)-\mv{B}\sin(\theta/2)\),
where \(\mv{B}\) is a unit bivector.

\item
In two dimensions, the even subalgebra \(a+b\mv{I}\)
behaves like the complex numbers.
Vectors are still represented as \(x \e1+y \e2\), not as \(x+y\mv{I}\).
The usual complex multiplication formula for rotations
is a simplification of the rotor sandwich product.

\item
In three dimensions, the even subalgebra is isomorphic to the quaternions.
Quaternion units correspond naturally to bivectors, or equivalently to axial
vectors after applying three-dimensional duality.

\item
In four and higher dimensions, general rotations need not be simple.
Their rotors may contain higher even-grade parts,
such as a grade-four part in four dimensions.
\end{itemize}

\section{Conclusion}
\label{sec:conclusion}
We compared four approaches to a concrete rotation problem in~\(\R^3\),
ranging from an ad hoc geometric argument to geometric algebra.
Although all methods produce the same result,
they emphasize different aspects of rotations.

The geometric-algebra approach requires the greatest conceptual investment,
but it also provides the broadest perspective.
Through rotors and the Cartan–Dieudonné theorem,
reflections, planar rotations, spatial rotations, quaternions, and
higher-dimensional rotations
all become manifestations of the same underlying structure.

\paragraph{Acknowledgments}
I used Mathematica for verification
(see Appendix~\ref{sec:mathematica-code})
and for creating the illustrations.
I had some helpful conversations with ChatGPT-4o and -o3-mini,
and for the 2026 update with ChatGPT-5.

\appendix

\section{Geometric Algebra details}
\label{sec:GA-details}

The geometric product of two vectors
\(\mv{v} = v_x \e{1} + v_y \e{2} + v_z \e{3}\)
and \(\mv{w} = w_x \e{1} + w_y \e{2} + w_z \e{3}\) is calculated as
\begin{align*}
  \mv{v} \mv{w} 
    ={} & \left(v_x \e1 + v_y \e2 + v_z \e3\right)
    \left(w_x \e1 + w_y \e2 + w_z \e3\right) \\
    ={} & v_x \e1 \left(w_x \e1 + w_y \e2 + w_z \e3\right) +{}\\
        & v_y \e2 \left(w_x \e1 + w_y \e2 + w_z \e3\right) +{}\\
        & v_z \e3 \left(w_x \e1 + w_y \e2 + w_z \e3\right) \\
    ={} & v_x w_x \e{11} + v_x w_y \e{12} + v_x w_z \e{13} +{}\\
        & v_y w_x \e{21} + v_y w_y \e{22} + v_y w_z \e{23} +{}\\
        & v_z w_x \e{31} + v_z w_y \e{32} + v_z w_z \e{33} \\
    ={} & \phantom{{}-{}}v_x w_x\phantom{\e{11}} + v_x w_y \e{12} - v_x w_z \e{31} +{}\\
        & -v_y w_x \e{12} + v_y w_y\phantom{\e{22}} + v_y w_z \e{23} +{}\\
        & \phantom{{}-{}}v_z w_x \e{31} - v_z w_y \e{23} + v_z w_z\phantom{\e{33}} \\
    ={} & v_x w_x + v_y w_y + v_z w_z +{}\\
        & \left(v_y w_z - v_z w_y\right) \e{23} +
          \left(v_z w_x - v_x w_z\right) \e{31} +
          \left(v_x w_y - v_y w_x\right) \e{12} \\
    ={} & \mv{v} \cdot \mv{w} + \mv{v} \wedge \mv{w}
\end{align*}
which defines \(\mv{v} \wedge \mv{w} = \mv{v}\mv{w} - \mv{v} \cdot \mv{w}\), called the \emph{wedge product}.
Note the resemblance and difference to \(\mv{v} \times \mv{w}\):
the coefficients are the same,
but \(\mv{v} \wedge \mv{w}\) uses the dual bivectors as basis
compared to the regular basis vectors used in \(\mv{v} \times \mv{w}\).
We do have
\[\|\mv{v} \wedge \mv{w}\| = \|\mv{v} \times \mv{w}\| =
  \|\mv{v}\| \|\mv{w}\| \sin \theta
\]
where \(\theta\) is the angle between~\(\mv{v}\) and~\(\mv{w}\).

Also note that in general the geometric product of two vectors
in three dimensions is the sum of a scalar and a bivector.
In fact, the product of two vectors 
is like a complex number
in the plane spanned by those vectors,
since such a product is a linear combination of the scalar unit~\(1\)
and a unit bivector,
which squares to~\(-1\).
To see the latter,
consider bivector \(x\, \e{23} + y\, \e{31} + z\, \e{12}\)
and calculate its square.
This comes to \(-\left(x^2+y^2+z^2\right)\).

The scalar part \(\mv{v} \cdot \mv{w}\) equals zero if and only if
\(\mv{v} \perp \mv{w}\)
and the bivector part \(\mv{v} \wedge \mv{w}\) vanishes if and only if
\(\mv{v} \parallel \mv{w}\).
We also see that if \(\mv{v} \perp \mv{w}\),
then \(\mv{v}\mv{w} = \mv{v} \wedge \mv{w} = - \mv{w} \wedge \mv{v} = - \mv{w}\mv{v}\).
Moreover, if \(\mv{v} \parallel \mv{w}\),
then \(\mv{v}\mv{w} = \mv{v} \cdot \mv{w} = \mv{w} \cdot \mv{v} = \mv{w}\mv{v}\).

\subsection{Formula for reflection}
\label{sec:GA-reflection}
The formula for reflecting vector~\(\mv{v}\) in
the plane perpendicular to unit vector~\(\mv{b}\)
(in general, in the hyperplane dual to~\(\mv{b}\))
was given as \(-\mv{b}\mv{v}\mv{b}\).
This can be explained as follows:
\begin{align*}
  -\mv{b}\mv{v}\mv{b} ={} & -\mv{b} \left(\mv{v}_{\perp \mv{b}} + \mv{v}_{\parallel \mv{b}}\right) \mv{b} &&
      \text{distribution}\\
    ={} & -\mv{b}\, \mv{v}_{\perp \mv{b}}\, \mv{b} - \mv{b}\, \mv{v}_{\parallel \mv{b}}\, \mv{b} &&
      \mv{v}_{\perp \mv{b}} \perp \mv{b} \text{ and } \mv{v}_{\parallel \mv{b}} \parallel \mv{b}\\
    ={} & \mv{b}^2\, \mv{v}_{\perp \mv{b}} - \mv{b}^2\, \mv{v}_{\parallel \mv{b}} &&
      \mv{b} \text{ is unit vector} \\
    ={} & \mv{v}_{\perp \mv{b}} - \mv{v}_{\parallel \mv{b}}
\end{align*}
So, indeed \(-\mv{b}\mv{v}\mv{b}\) flips the component of~\(\mv{v}\) parallel to~\(\mv{b}\).
(Similarly, \(\mv{b}\mv{v}\mv{b}\) flips the component of~\(\mv{v}\) perpendicular to~\(\mv{b}\).)

\subsection{Formula for rotor}
\label{sec:GA-rotor}
To understand the rotor \(R = e^{-\mv{B}\frac{\theta}{2}}\),
suppose we have unit vectors \(\mv{b}, \mv{c}\) spanning angle~\(\theta/2\).
Then we can calculate:
\begin{align*}
  \mv{c}\mv{b} ={} & \mv{c} \cdot \mv{b} + \mv{c} \wedge \mv{b} \\
    ={} & \mv{b} \cdot \mv{c} - \mv{b} \wedge \mv{c} \\
    ={} & \|\mv{b}\| \|\mv{c}\| \cos\tfrac{\theta}{2} -
      \|\mv{b}\| \|\mv{c}\| \frac{\mv{b} \wedge \mv{c}}{\|\mv{b} \wedge \mv{c}\|}\sin\tfrac{\theta}{2} \\
    ={} & \cos\tfrac{\theta}{2} - \mv{B} \sin\tfrac{\theta}{2}
      \quad\text{where } \mv{B} = \frac{\mv{b} \wedge \mv{c}}{\|\mv{b} \wedge \mv{c}\|} \\
    ={} & e^{- \mv{B} \frac{\theta}{2}}
\end{align*}
We express the bivector as \(\mv{b}\wedge\mv{c}\)
rather than \(\mv{c}\wedge\mv{b}\),
because this orientation is consistent with the chosen rotation axis
via the right-hand rule.
Also keep in mind that \(\mv{B}\) is a unit bivector
(hence, \(\mv{B}^2 = -1\))
in the plane spanned by~\(\mv{b}\) and~\(\mv{c}\),
which is used for rotating from~\(\mv{b}\) to~\(\mv{c}\).
To rotate in the plane spanned by~\(\mv{b}\) and~\(\mv{c}\)
by angle~\(\theta\),
any pair \(\mv{b}, \mv{c}\) spanning angle~\(\theta/2\) will do.
Thus, the unit bivector in the plane spanned by~\(\mv{b}\) and~\(\mv{c}\) will do,
as long as it is scaled to the right size by \(\sin\tfrac{\theta}{2}\).
In particular, we don't need concrete vectors~\(\mv{b}\) and~\(\mv{c}\).

\subsection{Formula for bivector perpendicular to vector}
\label{sec:GA-perpendicular-bivector}
When \(\mv{n}\) is a unit vector,
then \(\mv{I}\mv{n}\) is a unit bivector perpendicular to~\(\mv{n}\).
We can understand this as follows.
Choose the basis such that \(\mv{n}=\e1\).
Then we have \(\mv{I}\mv{n} = \e{123}\e1 = \e{1123} = \e{23}\),
which is a unit bivector that is spanned by~\(\e2\) and~\(\e3\),
both of which are perpendicular to~\(\e1=\mv{n}\).

\subsection{Full algebraic calculation}
\label{sec:GA-calculation}
Finally,
let's do the algebraic calculation alluded to earlier.
\begin{align*}
      R P \rev{R}
  ={} & \left(\tfrac{\sqrt{3}}{2} - \tfrac{1}{2}\, \mv{B}\right)
    (\e1 + \e3) \left(\tfrac{\sqrt{3}}{2} + \tfrac{1}{2}\, \mv{B}\right) \\
  ={} & (s - t(\e{23} + \e{31} + \e{12}))
    (\e1 + \e3) (s + t(\e{23} + \e{31} + \e{12}))\quad\text{where }
  s, t = \tfrac{\sqrt{3}}{2},\, \tfrac{1}{2\sqrt{3}} \\
  ={} & s^2(\e1 + \e3) +{} \\
      & st((\e1 + \e3)(\e{23} + \e{31} + \e{12}) -
           (\e{23} + \e{31} + \e{12})(\e1 + \e3)) +{}\\
      & t^2 (\e{23} + \e{31} + \e{12})(\e1 + \e3)(\e{23} + \e{31} + \e{12}) \\
  ={} & s^2(\e1 + \e3) +{} \\
      & st(\phantom{{}-{}}\e{123} + \e{131} + \e{112} + \e{323} + \e{331} + \e{312} +{}\\
      & \phantom{st(}{}- \e{231} - \e{311} - \e{121} - \e{233} - \e{313} - \e{123}) -{}\\
      & t^2 (\e{23}(\e1 + \e3)(\e{23} + \e{31} + \e{12}) +{}\\
      & \phantom{t^2 (}\e{31}(\e1 + \e3)(\e{23} + \e{31} + \e{12}) +{}\\
            & \phantom{t^2 (}\e{12}(\e1 + \e3)(\e{23} + \e{31} + \e{12})) \\
  ={} & s^2(\e1 + \e3) +{} \\
      & st(\e{123} - \e{3} + \e{2} - \e{2} + \e{1} + \e{123} -
           \e{123} - \e{3} + \e{2} - \e{2} + \e{1} - \e{123}) -{}\\
      & t^2 (\e{23123} + \e{23131} + \e{23112} + \e{23323} + \e{23331} + \e{23312} +{}\\
      & \phantom{t^2 (}\e{31123} + \e{31131} + \e{31112} + \e{31323} + \e{31331} + \e{31312} +{}\\
      & \phantom{t^2 (}\e{12123} + \e{12131} + \e{12112} + \e{12323} + \e{12331} + \e{12312}) \\
  ={} & s^2(\e1 + \e3) +
        2st(\e{1} - \e{3}) -{}\\
      & t^2 ({}-\e{1} - \e{2} - \e{3} + \e{3} + \e{123} - \e{1} +{}\\
      & \phantom{t^2 (}-\e{2} + \e{1} + \e{123} - \e{123} + \e{3} - \e{2} +{}\\
      & \phantom{t^2 (}-\e{3} - \e{123} + \e{1} - \e{1} - \e{2} - \e{3}) \\
  ={} & s^2(\e1 + \e3) +
        2st(\e{1} - \e{3}) +
        t^2 (\e{1} + 4\e{2} + \e{3}) \\
  ={} & (s^2 + 2st + t^2)\e1 + 4t^2 \e{2} + (s^2 - 2st + t^2)\e3 \\
  ={} & (s + t)^2 \e1 + 4t^2 \e{2} + (s - t)^2 \e3 \\
  ={} & \tfrac{4}{3} \e1 + \tfrac{1}{3} \e2 + \tfrac{1}{3} \e3
\end{align*}
This is a good exercise in algebra and sign bookkeeping,
which is ideal for automation.
    
\section{Mathematica code}
\label{sec:mathematica-code}

The general set up of our rotation problem:
\begin{lstlisting}
origin = {0, 0, 0}
axis = {1, 1, 1}
theta = 60 Degrees (* just an angle; we should also do it for an arbitrary angle *)
pointP = {1, 0, 1}
\end{lstlisting}

Method~1:
\begin{lstlisting}
pointP2 = {1, 1, 0};
pointP3 = {0, 1, 1};
pointQ = Mean[{pointP, pointP2, pointP3}];
pointM = Mean[{pointP, pointP2}];
pointPdash = pointQ + (pointM - pointQ) Norm[pointP - pointQ] / Norm[pointM - pointQ]
\end{lstlisting}

Method~2:
\begin{lstlisting}
zAxis = Normalize[axis];
pointQ = (pointP . zAxis) zAxis; (* using dot product of vectors *)
xAxis = Normalize[pointP - pointQ];
yAxis = Cross[zAxis, xAxis];
pointPdash = pointQ + (Cos[theta] xAxis + Sin[theta] yAxis) Norm[pointP - pointQ]
\end{lstlisting}

Method~3:
\begin{lstlisting}
R = RotationMatrix[theta, axis];
pointPdash = R . pointP (* using matrix-vector product *)
\end{lstlisting}

Method~4:
\begin{lstlisting}
<< CliffordBasic` (* installed from https://github.com/jlaragonvera/Geometric-Algebra *)
a = ToBasis[axis] (* axis: e[1] + e[2] + e[3] *)
p = ToBasis[pointP] (* point to be rotated: e[1] + e[3] *)
B = -Dual[a, 3] (* bivector for rotation plane: e[2,3] + e[3,1] + e[1,2] *)
rotor[B_, theta_] := Cos[theta/2] - B Sin[theta/2] / Magnitude[B];
GeometricRotate[v_, R_] := GeometricProduct[R, v, Turn[R]]; (* Turn[R] = reverse of R *)
pdash = GeometricRotate[p, rotor[B, theta]];
pointPdash = ToVector[pdash]
\end{lstlisting}

\end{document}